\documentclass[11pt]{amsart}

\usepackage{amsmath,amsthm,amsfonts,amssymb,tikz,extarrows}
\usepackage{hyperref}
\hypersetup{linkbordercolor={1 .8 .8},citebordercolor={.7 1 .7}}
\newcommand{\prob}{\mathbb{P}}

\newcommand{\real}{\mathbb{R}}
\newcommand{\C}{\mathbb{C}}
\newcommand{\N}{\mathbb{N}}

\newcommand{\Inte}{\mathbb{Z}}
\newcommand{\ind}{\mathbb{I}}
\newcommand{\ve}{e}

\newcommand{\sh}{\text{sh}}
\newcommand{\ty}{\text{ty}}
\newcommand{\gl}{\mathfrak{gl}}

\newcommand{\vars}{\mathcal S}
\newcommand{\vart}{\mathcal T}

\newtheorem{thm}{Theorem}
\newtheorem{prop}[thm]{Proposition}

\newtheorem{cly}[thm]{Corollary}

\theoremstyle{definition}
\newtheorem{exa}{Example}

\begin{document}

\title[A $q$-weighted version of the Robinson-Schensted algorithm]{A $q$-weighted version of the Robinson-Schensted algorithm}

\author{Neil O'Connell}
\address{Mathematics Institute, University of Warwick, Coventry CV4 7AL, UK}
\email{n.m.o-connell@warwick.ac.uk}
\author{Yuchen Pei}
\address{Mathematics Institute, University of Warwick, Coventry CV4 7AL, UK}
\email{y.pei@warwick.ac.uk}




\maketitle

\begin{abstract}
We introduce a $q$-weighted version of the Robinson-Schensted (column insertion) algorithm 
which is closely connected to $q$-Whittaker functions (or Macdonald polynomials with $t=0$) and 
reduces to the usual Robinson-Schensted algorithm when $q=0$.  The $q$-insertion algorithm is `randomised', 
or `quantum', in the sense that when inserting a 
positive integer into a tableau, the output is a distribution of weights on a particular set of tableaux which 
includes the output which would have been obtained via the usual column insertion algorithm.   
There is also a notion of recording tableau in this setting.  
We show that the distribution of weights of the pair of tableaux obtained when 
one applies the $q$-insertion algorithm to a random word or permutation takes a particularly 
simple form and is closely related to $q$-Whittaker functions.
In the case $0\le q<1$, the $q$-insertion algorithm applied to a random word also provides a new
framework for solving the $q$-TASEP interacting particle system introduced (in the language of $q$-bosons) 
by Sasamoto and Wadati \cite{sw} and yields formulas which are equivalent
to some of those recently obtained by Borodin and Corwin \cite{bc} via a stochastic evolution on 
discrete Gelfand-Tsetlin patterns (or semistandard tableaux) which is coupled to the $q$-TASEP process.  
We show that the sequence of $P$-tableaux obtained when one applies the $q$-insertion algorithm to a random 
word defines another, quite different, evolution on semistandard tableaux which is also coupled to the $q$-TASEP 
process.
\end{abstract}


\section{Introduction}

We introduce a $q$-weighted version of the Robinson-Schensted (column insertion) algorithm 
which is closely connected to $q$-Whittaker functions 
(or Macdonald polynomials with $t=0$) and reduces to the usual Robinson-Schensted algorithm 
when $q=0$.   The insertion algorithm is `randomised', or `quantum', in the sense 
that when inserting a positive integer into a tableau, the output is a distribution of weights 
on a particular set of tableau which includes the output which would have been obtained
via the usual column insertion algorithm.  As such, it is similar to the quantum insertion
algorithm introduced by Date, Jimbo and Miwa \cite{djm} (see also \cite{berg}) but with different weights.  
There is also a notion of recording tableau in this setting.  
We show that the distribution of weights of the pair of tableaux obtained when 
one applies the insertion algorithm to a random word or permutation takes a particularly 
simple form and is closely related to $q$-Whittaker functions. These are functions
defined on integer partitions which are eigenfunctions the relativistic Toda chain~\cite{ruijsenaars,ruijsenaars1,etingof,glo1} 
and simply related to Macdonald polynomials (as a function of the index) with the parameter $t=0$~\cite{glo3}.
When $q=0$, they are given by Schur polynomials.  Our main result provides a starting point for developing a 
new combinatorial framework for $q$-Whittaker functions and related objects, such as Demazure and 
Kirillov-Reshetikhin crystals.  It will be interesting to understand the relation to recent 
developments in this area, see~\cite{hhl,l,ry,bbl,ll,st,bf,ls} and references therein.

In the case $0\le q<1$, the $q$-insertion algorithm applied to a random word also provides a new
framework for solving the $q$-TASEP interacting particle system introduced (in the language of $q$-bosons) 
by Sasamoto and Wadati \cite{sw} and yields formulas which are equivalent
to some of those recently obtained by Borodin and Corwin \cite{bc} via a stochastic evolution on 
discrete Gelfand-Tsetlin patterns---or, equivalently, semistandard tableaux---which is coupled to the $q$-TASEP process.  
We show that the sequence of $P$-tableaux obtained when one applies the $q$-insertion algorithm to a random 
word defines another, quite different, evolution on semistandard tableaux which is also coupled to the $q$-TASEP 
process (after Poissonisation).  The $q$-TASEP process is a particular case of the totally asymmetric 
zero-range process~\cite{bks}.  See also~\cite{bcs} for related recent work.

When $q\to1$, the $q$-Whittaker functions converge with appropriate rescaling to $\gl_l$-Whittaker 
functions \cite{glo2}.   The main result of the present paper can be regarded as a natural (yet non-obvious)
discretisation, in time and space, of the main result of the paper~\cite{noc12}, 
which relates a continuous-time version of the geometric Robinson-Schensted-Knuth (RSK) 
correspondence introduced by 
A.N. Kirillov~\cite{kirillov1}, with Brownian motion as input, to the open quantum Toda chain with $l$ particles.  
A discrete time version of that result has been developed in the papers~\cite{cosz,osz},
which is formulated directly in the context of Kirillov's geometric RSK correspondence. 
The present work differs significantly from~\cite{noc12,cosz,osz} in that the analogue of the RSK mapping we
consider here is (necessarily) randomised.  In the above scaling limit, the $q$-insertion algorithm we introduce 
in this paper should converge in an appropriate sense to the continuous-time version of the geometric RSK 
mapping considered in~\cite{noc12}, which is deterministic, and the main result of this paper should rescale to 
the main result of~\cite{noc12}.  This can be seen by comparing with the corresponding scaling limits considered 
in~\cite{bc,glo2}.

The outline of the paper is as follows. In the next section we give some background on the
Robinson-Schensted algorithm.  In Section 3, we describe the $q$-weighted version of this algorithm.  
The main result is presented in Section 4.
In Section 5, we consider the $q$-insertion algorithm with $0\le q<1$ applied to a random word
and explain the connection to the $q$-TASEP interacting particle system.  In Section 6 we consider 
the algorithm applied to a random permutation.  The proofs are given in Section 7.

\section{The Robinson-Schensted algorithm}

The Robinson-Schensted algorithm is a combinatorial algorithm which plays
a fundamental role in the theory of Young tableaux~\cite{robinson,schensted,fulton,sagan,stanley}.  
There are two versions, which are in some sense dual to each other, defined via
insertion (or `bumping') algorithms known as {\em row insertion} and {\em column insertion}.  
The column insertion algorithm is also sometimes referred to as the dual RSK algorithm,
because it has a natural extension to zero-one matrices which was introduced by Knuth~\cite{knuth}.
It is the column insertion version which we consider and generalise in this paper.

A tableau $P$ is a Young diagram with positive integer entries which are weakly
increasing in each row and strictly increasing in each column.  The corresponding diagram
represents an integer partition which is referred to as the shape of the tableau $P$ and
denoted by $\sh P$. For example,
$$\begin{array}{ccccc}1&1&2&3&\\2&3&3&&\\3&&&&\end{array}$$
is a tableau with shape $(4,3,1)$.
To insert a positive integer $k$ into a tableau $P$, we begin by trying to place
that integer at the bottom of the first column of $P$.  If the result is a tableau, we are done.
Otherwise, it bumps the smallest entry in that column which is larger than or equal to $k$.
Now proceed by inserting the bumped entry into the second column according to the same rule, 
and so on, until we have placed a bumped entry at the bottom of column (or on its own in a new column).  
For example, if we insert the number 2 into the tableau shown above, the outcome is
$$\begin{array}{ccccc}1&1&2&3&3\\2&2&3&&\\3&&&&\end{array}$$
In this example, the 2 in the first column is bumped into the second, 
the 2 in the second is bumped into the third,
the 3 in the third column is bumped into the fourth,
and the 3 in the fourth is bumped into a new fifth column on its own.
Actually, it will be helpful for later reference to summarise this sequence of events in the following way:  
in this example, a 2 is inserted into the second row, and a 3 is bumped from the second row and inserted into the first row.

Now, applying this insertion algorithm recursively to a word $w=w_1\ldots w_n\in[l]^n$,
starting with an empty tableau and successively inserting the numbers $w_1,w_2,\ldots,w_n$,
gives rise to a sequence of tableau $P(1),P(2),\ldots,P(n)=P$.  Note that it is not possible in general
to recover the word $w$ from the tableau $P$.  This motivates the notion of a {\em recording tableau},
which we denote by $Q$.  The tableau $Q$ has size $n$ and is {\em standard}, that is, it contains each
of the numbers $1,2,\ldots,n$ exactly once.  If we denote by $Q^i$ the sub-tableau of $Q$ consisting only 
of those entries which are not greated than $i$, then $Q$ is defined by the requirement that 
$\sh Q^i=\sh P(i)$ for $1\le i\le n$.  
For example, if $w=1143232$ then
  \begin{align*}
  P=\begin{array}{cccc}
    1&1&3&4\\
    2&2&&\\
    3&&&
  \end{array}\qquad
  Q=\begin{array}{cccc}
    1&2&5&7\\
    3&4&&\\
    6&&&
  \end{array}
\end{align*}
The mapping $w\mapsto (P,Q)$ defines a bijection from the set of words $[l]^n$ to the set of pairs 
$(P,Q)\in\vart_l\times\vars_n$ such that $\sh P=\sh Q$, where $\vart_l$ denotes the set of tableaux 
with entries from $[l]$ and $\vars_n$ denotes the set of standard tableaux of size $n$.  It is the column
insertion version of the {\em Robinson-Schensted correspondence}.  

\begin{figure}[h!]
\begin{center}
\begin{tikzpicture}
  \foreach \i in {0,...,4}{
  \path(\i,0) node[below]{\i};
  }
  \foreach \j in {0,...,2}{
  \path(0,\j) node[left]{\j};
  }
  \draw[step=1cm,color=gray] (0,0) grid (4,2);
  \draw[->,line width=2pt] (2,2) -- (3,2) -- (3,1) -- (4,1);
\end{tikzpicture}
\end{center}
\caption{An insertion path}
\label{f:columnlattice}
\end{figure}
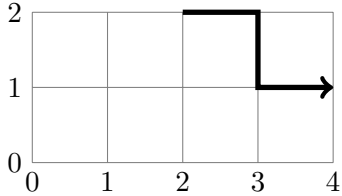

As a warm up for next section, we note that the above column insertion algorithm can also be described 
in terms of lattice paths, as follows.
Suppose we are inserting a number $k$ with $1\le k\le l$ into a tableau $P\in\mathcal T_l$, with resulting tableau $\tilde P$. 
For $1\le i\le l$, set $\lambda^i=\sh P^i$, and $\tilde\lambda^i=\sh\tilde P^i$.
Let $(e_i,\ 1\le i\le l)$ denote the standard basis in $\mathbb Z^l$. 
Then $\tilde\lambda^i=\lambda^i+e_{j_i}$ where $k=j_{k-1}\ge j_k\ge \dots \ge j_l\ge1$ is a weakly decreasing
sequence defined by
\begin{align*}
  j_i=\max\{\{2\le m\le j_{i-1}:\lambda^{i-1}_{m-1}-\lambda^i_m>0\}\cup\{1\}\},\quad i=k,k+1,\dots,l.
\end{align*}
The sequence $k\ge j_k\ge j_{k+1}\ge \cdots j_l\ge 1$ determines a down/right lattice path in $\Inte^2$ from 
$(k,k)$ to $(l+1,j_l)$ by specifying the $y$-coordinates at which the path moves to the right.  From the definition,
this path takes a horizontal step to the right $(i,j)\to(i+1,j)$ whenever $\lambda^{i-1}_{j-1}>\lambda^i_j$ 
or $j=1$, otherwise it takes a step down $(i,j)\to(i,j-1)$.  We will refer to this lattice path as the {\em insertion path}.
The interpretation is as follows.  A horizontal portion of the path starting at $(i,j)$ represents inserting an
$i$ into the $j$th row.  A vertical portion starting at $(i,j)$ and ending at $(i,j-r)$ indicates that an
$i$ is bumped from the $j$th row to the $(j-r)$th row.  For example, the insertion path corresponding to 
the previous example of inserting a 2 into the tableau
$$\begin{array}{ccccc}1&1&2&3&\\2&3&3&&\\3&&&&\end{array}$$
with $l=3$ is illustrated in Figure \ref{f:columnlattice}.

\section{The $q$-weighted version}

In this paper, we consider the following generalisation of the column insertion algorithm.  
It is defined by a collection of kernels $I_k(P,\tilde P)$ which depend on a complex parameter $q$.
We assume throughout that $q$ is not a root of unity.  
If $0\le q<1$, we interpret the quantity $I_k(P,\tilde P)$ as the probability that, 
when we insert $k$ into the tableau $P$, the output is $\tilde P$.  Recall that the {\em type} of a tableaux $P$, 
which we denote $\ty P$, is the composition $\mu=(\mu_1,\mu_2,\ldots)$ where $\mu_i$ is the number of $i$'s in $P$.
The set of $\tilde P$ for which $I_k(P,\tilde P)\ne 0$ has the following properties. The type of $\tilde P$
is given by $\ty \tilde P=\ty P+e_k$.  The shape of $\tilde P$ satisfies $\sh \tilde P=\sh P+e_j$ for some
$1\le j\le k$.  Moreover, if we set $\lambda^i=\sh P^i$ and $\tilde \lambda^i=\sh \tilde P^i$, then there is
a weakly decreasing sequence $k=j_{k-1}\ge j_k\ge j_{k+1}\ge \cdots j_l\ge 1$ such that 
$\tilde \lambda^i = \lambda^i$ for $1\le i<k$ and $\tilde \lambda^i = \lambda^i+e_{j_i}$ for $k\le i\le l$.
The kernel $I_k(P,\tilde P)$ is defined to be zero if there is no such sequence; if there is such a sequence,
it is given as follows.
Define 
\begin{align*}
f_0(i,j)=1-q^{\lambda^{i-1}_{j-1}-\lambda^i_j},\qquad 
f_1(i,j)=\frac{1-q^{\lambda^{i-1}_{j-1}-\lambda^i_j}}{1-q^{\lambda^{i-1}_{j-1}-\lambda^{i-1}_j}},\quad\text{ for }j>1;\\
f_0(i,1)=f_1(i,1)=1.
\end{align*}
and set
\begin{align*}
  f(i,j)=\begin{cases}
    f_1(i,j),&\text{ if }j=j_{i-1}\text{ and }i\neq k;\\
    f_0(i,j),&\text{ otherwise.}
  \end{cases}
\end{align*}
Then
\begin{equation}\label{Ik}
I_k(P,\tilde P) = \prod_{i=k}^l \left( f(i,j_i) \prod_{j=j_i+1}^{j_{i-1}} (1-f(i,j))\right) .
\end{equation}
It follows easily from the definition that
$$\sum_{\tilde P} I_k(P,\tilde P)=1.$$
If $0\le q<1$, then $I_k(P,\tilde P)\ge 0$.  In this case, for each $k$ and $P$, $I_k(P,\cdot)$ defines a probability
distribution on $\vart_l$ and we interpret $I_k(P,\tilde P)$ as the probability that, 
when we insert $k$ into the tableau $P$, the output is $\tilde P$.

The formula \eqref{Ik} can be interpreted in terms of insertion paths, as follows.
The sequence $k\ge j_k\ge j_{k+1}\ge \cdots j_l\ge 1$ determines a down/right lattice path
in $\Inte^2$ from $(k,k)$ to the vertical boundary $\{(l+1,j),\ 1\le j\le k\}$ by specifying 
the $y$-coordinates at which the path moves to the right.  The edge weights are $f(i,j)$ on the horizontal edge 
$(i,j)\to (i+1,j)$ and $1-f(i,j)$ on the vertical edge $(i,j)\to (i,j-1)$, and taking a product of these weights along the
path gives the weight $I_k(P,\tilde P)$ for the corresponding output $\tilde P$.
We interpret this path as the insertion path associated with $q$-inserting the number $k$ into $P$ with 
resulting tableau $\tilde P$.  As before, a horizontal portion of the path starting at $(i,j)$ represents inserting an
$i$ into the $j$th row.  A vertical portion starting at $(i,j)$ and ending at $(i,j-r)$ indicates that an
$i$ is bumped from the $j$th row to the $(j-r)$th row. 
When $q=0$, there is only 
one output tableau $\tilde P$ with non-zero weight, namely the output of the usual column insertion algorithm. 
Moreover, if we denote by $\omega_0$ the insertion path corresponding to this tableau and by $S(k,P)$ the set of insertion paths corresponding to the support of $I_k(P,\cdot)$ for nonzero $q$, then $\omega_0\in S(k,P)$ and it is the `highest' path in 
$S(k,P)$ in the sense that the sequence $k\ge j_k\ge j_{k+1}\ge \cdots j_l\ge 1$ is maximal (in the second example below,
it is the path shown on the top left of Figure 2).

Let us compute the kernel $I_k(P,\tilde P)$ for some concrete examples.
\begin{exa}\label{ex:1} Suppose $l=2$.  If we are inserting a 1 into $P\in\vart_2$ there is only one possible outcome $\tilde P$
with $I_1(P,\tilde P)\ne 0$,
namely the one obtained by the usual column insertion algorithm: the 1 is inserted into the first row, pushing
the existing first row over by one. The weighted insertion path in this case is very simple:
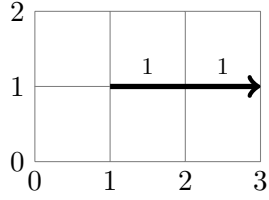
\[
\begin{tikzpicture}
  \foreach \i in {0,...,3}{
  \path(\i,0) node[below]{\i};
  }
  \foreach \j in {0,...,2}{
  \path(0,\j) node[left]{\j};
  }
  \draw[step=1cm,color=gray] (0,0) grid (3,2);
  \draw[->,line width=2pt,font=\footnotesize] (1,1) -- node[above]{1} (2,1) -- node[above]{1} (3,1);
\end{tikzpicture}
\]
For example, if
$$P=\begin{array}{cccc}1&1&2&2\\2&&&\end{array}$$
then, setting
$$\tilde P_1=\begin{array}{ccccc}1&1&1&2&2\\2&&&&\end{array}$$
we have
$$I_1(P,\tilde P)=\begin{cases} 1 & \mbox{ if } \tilde P=\tilde P_1\\0& \mbox{ otherwise.}\end{cases}$$
On the other hand, if we are inserting a 2 there are two possibilities:  \begin{enumerate}\item The 2 is inserted into the
second row, pushing the existing second row over by one: this outcome has weight $1-q^{\lambda^1_1-\lambda^2_2}$.
\item The 2 is inserted into the first row, pushing the existing 2's over by one:  this outcome has weight $q^{\lambda^1_1-\lambda^2_2}$.
\end{enumerate}
Note that these weights sum to one, as is always the case.
The corresponding insertion paths, with edge weights indicated, are:
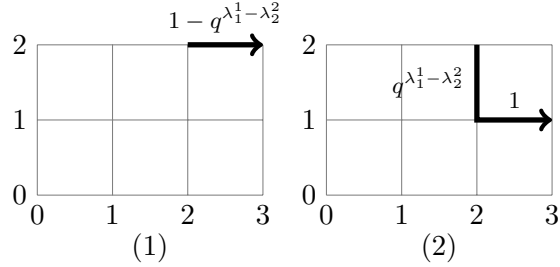
\[
\begin{tikzpicture}
  \foreach \i in {0,...,3}{
  \path(\i,0) node[below]{\i};
  }
  \foreach \j in {0,...,2}{
  \path(0,\j) node[left]{\j};
  }
  \draw[step=1cm,color=gray] (0,0) grid (3,2);
  \draw[->,line width=2pt,font=\footnotesize] (2,2) -- node[above]{$1-q^{\lambda^1_1-\lambda^2_2}$} (3,2); 
  \path(1.5,-0.7) node{(1)};
\end{tikzpicture}
\begin{tikzpicture}
  \foreach \i in {0,...,3}{
  \path(\i,0) node[below]{\i};
  }
  \foreach \j in {0,...,2}{
  \path(0,\j) node[left]{\j};
  }
  \draw[step=1cm,color=gray] (0,0) grid (3,2);
  \draw[->,line width=2pt,font=\footnotesize] (2,2) -- node[left]{$q^{\lambda^1_1-\lambda^2_2}$} (2,1) -- node[above]{1} (3,1); 
  \path(1.5,-0.7) node{(2)};
\end{tikzpicture}
\]
The quantity $\lambda^1_1-\lambda^2_2$ is the difference between the number of 1's in the first row and the number of 2's in the
second row, see Figure \ref{fig:ex1}.  
\begin{figure}[h!]
  \begin{center}
    \begin{tikzpicture}[font=\footnotesize,scale=.7,>=stealth]
      \draw (0,0) rectangle node{1} (1,-1);
      \draw (1,0) rectangle node{\dots} (6,-1);
      \draw (6,0) rectangle node{1} (7,-1);
      \draw (7,0) rectangle node{2} (8,-1);
      \draw (8,0) rectangle node{\dots} (10,-1);
      \draw (10,0) rectangle node{2} (11,-1);
      
      \draw (0,-1) rectangle node{2} (1,-2);
      \draw (1,-1) rectangle node{\dots} (3,-2);
      \draw (3,-1) rectangle node{2} (4,-2);
      \draw [dashed] (4,-2) -- (4,-2.5) (7,-1) -- (7,-2.5);
      \draw [<->] (4,-2) -- node[above]{$\lambda^1_1-\lambda^2_2$} (7,-2);
    \end{tikzpicture}
  \end{center}
  \caption{The quantity $\lambda^1_1-\lambda^2_2$ in the exponent in Example \ref{ex:1}}
  \label{fig:ex1}
\end{figure}
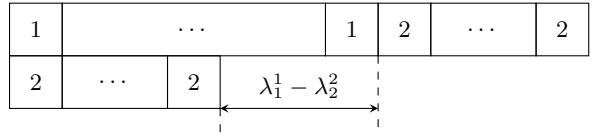

For example, inserting a 2 into 
$$P=\begin{array}{cccc}1&1&2&2\\2&&&\end{array}$$
gives
$$I_2(P,\tilde P)=\begin{cases} 1-q &  \mbox{ if } \tilde P=\tilde P_2\\ q &  \mbox{ if } \tilde P=\tilde P_3 \\
0& \mbox{ otherwise.} \end{cases}$$
where
$$\tilde P_2=\begin{array}{cccc}1&1&2&2\\2&2&&\end{array}$$
and
$$\tilde P_3=\begin{array}{ccccc}1&1&2&2&2\\2&&&&\end{array}.$$
\end{exa}
\begin{exa}\label{ex:2} Suppose $l=3$.  
   If we are inserting a 1 into $P\in\vart_3$ there is only one possible outcome $\tilde P$
with $I_1(P,\tilde P)\ne 0$,
namely the one obtained by the usual column insertion algorithm: the 1 is inserted into the first row, pushing
the existing first row over by one. The corresponding weighted insertion path is:
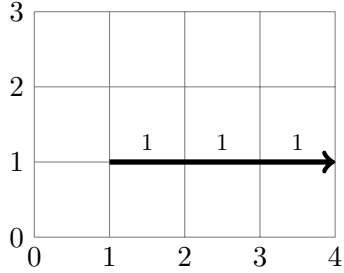
\[
\begin{tikzpicture}
  \foreach \i in {0,...,4}{
  \path(\i,0) node[below]{\i};
  }
  \foreach \j in {0,...,3}{
  \path(0,\j) node[left]{\j};
  }
  \draw[step=1cm,color=gray] (0,0) grid (4,3);
  \draw[->,line width=2pt,font=\footnotesize] (1,1) -- node[above]{1} (2,1) -- node[above]{1} (3,1)-- node[above]{1} (4,1); 
\end{tikzpicture}
\]

If we are inserting a 2, there are three possible outcomes:
\begin{enumerate}
\item The 2 is inserted into the second row, pushing the existing second row over by one: this outcome has weight 
$$(1-q^{\lambda^1_1-\lambda^2_2}) \dfrac{1-q^{\lambda^2_1-\lambda^3_2}}{1-q^{\lambda^2_1-\lambda^2_2}} ;$$
\item The 2 is inserted into the second row, bumping a 3 into the first row: 
 this outcome has weight 
$$(1-q^{\lambda^1_1-\lambda^2_2}) \left( 1-\dfrac{1-q^{\lambda^2_1-\lambda^3_2}}{1-q^{\lambda^2_1-\lambda^2_2}}\right) ;$$
\item The 2 is inserted into the first row, pushing existing 2's and 3's in first row over by one: this outcome has weight 
$q^{\lambda^1_1-\lambda^2_2}$.
\end{enumerate}
The corresponding insertion paths, with edge weights indicated, are:
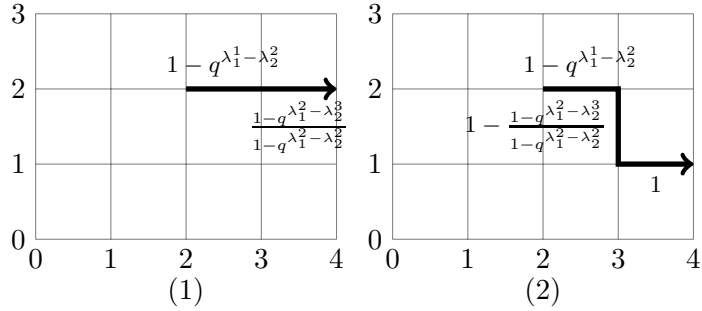
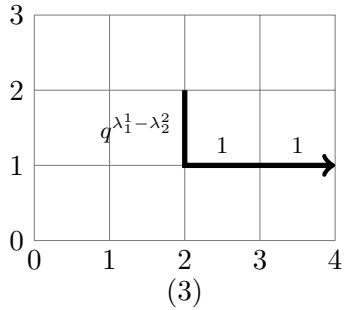
\[
\begin{tikzpicture}
  \foreach \i in {0,...,4}{
  \path(\i,0) node[below]{\i};
  }
  \foreach \j in {0,...,3}{
  \path(0,\j) node[left]{\j};
  }
  \draw[step=1cm,color=gray] (0,0) grid (4,3);
  \draw[->,line width=2pt,font=\footnotesize] (2,2) -- node[above]{$1-q^{\lambda^1_1-\lambda^2_2}$} (3,2) -- node[below]{$\frac{1-q^{\lambda^2_1-\lambda^3_2}}{1-q^{\lambda^2_1-\lambda^2_2}}$} (4,2); 
  \path(2,-0.7) node{(1)};
\end{tikzpicture}
\begin{tikzpicture}
  \foreach \i in {0,...,4}{
  \path(\i,0) node[below]{\i};
  }
  \foreach \j in {0,...,3}{
  \path(0,\j) node[left]{\j};
  }
  \draw[step=1cm,color=gray] (0,0) grid (4,3);
  \draw[->,line width=2pt,font=\footnotesize] (2,2) -- node[above]{$1-q^{\lambda^1_1-\lambda^2_2}$} (3,2) -- node[left]{$1-\frac{1-q^{\lambda^2_1-\lambda^3_2}}{1-q^{\lambda^2_1-\lambda^2_2}}$} (3,1) -- node[below]{1} (4,1); 
  \path(2,-0.7) node{(2)};
\end{tikzpicture}
\]
\[
\begin{tikzpicture}
  \foreach \i in {0,...,4}{
  \path(\i,0) node[below]{\i};
  }
  \foreach \j in {0,...,3}{
  \path(0,\j) node[left]{\j};
  }
  \draw[step=1cm,color=gray] (0,0) grid (4,3);
  \draw[->,line width=2pt,font=\footnotesize] (2,2) -- node[left]{$q^{\lambda^1_1-\lambda^2_2}$} (2,1) -- node[above]{$1$} (3,1) -- node[above]{1} (4,1); 
  \path(2,-0.7) node{(3)};
\end{tikzpicture}
\]

If we are inserting a 3, there are also three possible outcomes: the 3 is placed in the third, second or first row with
respective weights $1-q^{\lambda^2_2-\lambda^3_3}$, $q^{\lambda^2_2-\lambda^3_3}(1-q^{\lambda^2_1-\lambda^3_2})$
and $q^{\lambda^2_2-\lambda^3_3}q^{\lambda^2_1-\lambda^3_2}$. The corresponding insertion paths, with edge weights
indicated, are: 
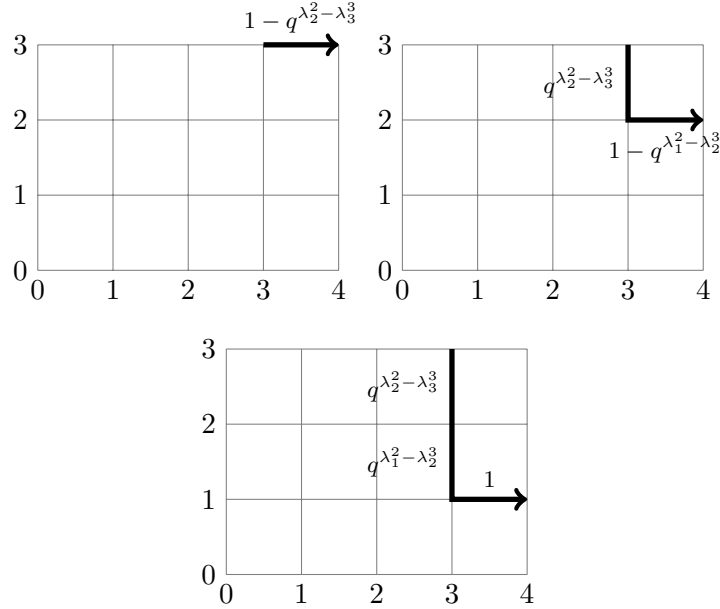
\[
\begin{tikzpicture}
  \foreach \i in {0,...,4}{
  \path(\i,0) node[below]{\i};
  }
  \foreach \j in {0,...,3}{
  \path(0,\j) node[left]{\j};
  }
  \draw[step=1cm,color=gray] (0,0) grid (4,3);
  \draw[->,line width=2pt,font=\footnotesize] (3,3) -- node[above]{$1-q^{\lambda^2_2-\lambda^3_3}$} (4,3); 
\end{tikzpicture}
\begin{tikzpicture}
  \foreach \i in {0,...,4}{
  \path(\i,0) node[below]{\i};
  }
  \foreach \j in {0,...,3}{
  \path(0,\j) node[left]{\j};
  }
  \draw[step=1cm,color=gray] (0,0) grid (4,3);
  \draw[->,line width=2pt,font=\footnotesize] (3,3) -- node[left]{$q^{\lambda^2_2-\lambda^3_3}$} (3,2) -- node[below]{$1-q^{\lambda^2_1-\lambda^3_2}$} (4,2); 
\end{tikzpicture}
\]
\[
\begin{tikzpicture}
  \foreach \i in {0,...,4}{
  \path(\i,0) node[below]{\i};
  }
  \foreach \j in {0,...,3}{
  \path(0,\j) node[left]{\j};
  }
  \draw[step=1cm,color=gray] (0,0) grid (4,3);
  \draw[->,line width=2pt,font=\footnotesize] (3,3) -- node[left]{$q^{\lambda^2_2-\lambda^3_3}$} (3,2) -- node[left]{$q^{\lambda^2_1-\lambda^3_2}$} (3,1) -- node[above]{1} (4,1); 
\end{tikzpicture}
\]

The quantities $\lambda^1_1-\lambda^2_2$, $\lambda^2_1-\lambda^3_2$, etc. which appear in the above weights
are illustrated in Figure \ref{fig:ex2}.
\begin{figure}[h!]
  \begin{center}
    \begin{tikzpicture}[font=\footnotesize,scale=.7,>=stealth]
      \draw (0,0) rectangle node{1} (1,-1);
      \draw (1,0) rectangle node{\dots} (7,-1);
      \draw (7,0) rectangle node{1} (8,-1);
      \draw (8,0) rectangle node{2} (9,-1);
      \draw (9,0) rectangle node{\dots} (12,-1);
      \draw (12,0) rectangle node{2} (13,-1);
      \draw (13,0) rectangle node{3} (14,-1);
      \draw (14,0) rectangle node{\dots} (16,-1);
      \draw (16,0) rectangle node{3} (17,-1);
      
      \draw (0,-1) rectangle node{2} (1,-2);
      \draw (1,-1) rectangle node{\dots} (4,-2);
      \draw (4,-1) rectangle node{2} (5,-2);
      \draw (5,-1) rectangle node{3} (6,-2);
      \draw (6,-1) rectangle node{\dots} (10,-2);
      \draw (10,-1) rectangle node{3} (11,-2);

      \draw (0,-2) rectangle node{3} (1,-3);
      \draw (1,-2) rectangle node{\dots} (2,-3);
      \draw (2,-2) rectangle node{3} (3,-3);
      \draw [dashed] (5,-2) -- (5,-3.5) (8,-1) -- (8,-3.5);
      \draw [<->] (5,-3) -- node[above]{$\lambda^1_1-\lambda^2_2$} (8,-3);
      \draw [dashed] (11,-2) -- (11,-3.5) (13,-1) -- (13,-3.5);
      \draw [<->] (11,-3) -- node[above]{$\lambda^2_1-\lambda^3_2$} (13,-3);
      \draw [dashed] (3,-3) -- (3,-4.5) (5,-2) -- (5,-4.5);
      \draw [<->] (3,-4) -- node[above]{$\lambda^2_2-\lambda^3_3$} (5,-4);
      \draw [dashed] (5,-2) -- (5,-4.5) (13,-2) -- (13,-4.5);
      \draw [<->] (5,-4) -- node[above]{$\lambda^2_1-\lambda^2_2$} (13,-4);
    \end{tikzpicture}
  \end{center}
  \caption{The exponent quantities in Example \ref{ex:2}.}
  \label{fig:ex2}
\end{figure}
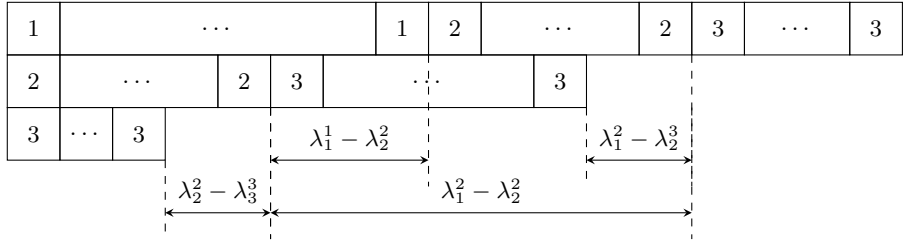
\end{exa}
\begin{exa} Suppose we are inserting a $3$ into 
\begin{align}\label{P}
  P=
  \begin{array}{cccccc}
    1&2&2&2&3&5\\
    2&3&4&5&&\\
    3&4&&&&\\
    5&&&&&
  \end{array}
\end{align}
The (four) possible output tableaux $\tilde P$ and their weights $I_3(P,\tilde P)$ are
shown in Figure \ref{f:qlattice}, along with the corresponding weighted insertion paths.
\end{exa}

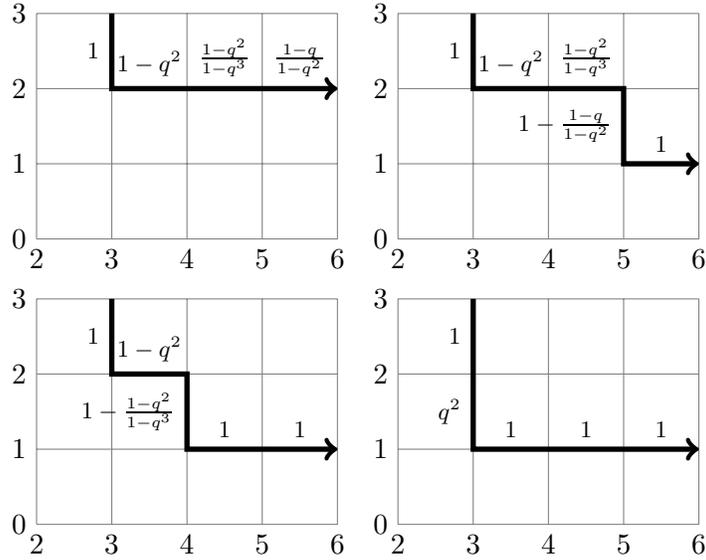
\begin{figure}[h!]

\begin{center}

$$  \tilde P=
    \begin{array}{cccccc}
      1&2&2&2&3&5\\
      2&3&3&4&5&\\
      3&4&&&&\\
      5&&&&&
    \end{array},\qquad I_3(P,\tilde P) = (1-q^2)\frac{1-q^2}{1-q^3}\frac{1-q}{1-q^2};$$
$$  \tilde P=
    \begin{array}{ccccccc}
      1&2&2&2&3&5&5\\
      2&3&3&4&&&\\
      3&4&&&&&\\
      5&&&&&&
    \end{array},\qquad I_3(P,\tilde P) = (1-q^2)\frac{1-q^2}{1-q^3}\left(1-\frac{1-q}{1-q^2}\right);$$
$$  \tilde P=    \begin{array}{ccccccc}
      1&2&2&2&3&4&5\\
      2&3&3&5&&&\\
      3&4&&&&&\\
      5&&&&&&
    \end{array},\qquad I_3(P,\tilde P) = (1-q^2)\left(1-\frac{1-q^2}{1-q^3}\right);$$
$$  \tilde P=    
    \begin{array}{ccccccc}
      1&2&2&2&3&3&5\\
      2&3&4&5&&&\\
      3&4&&&&&\\
      5&&&&&&
    \end{array},\qquad I_3(P,\tilde P) = q^2.$$

\begin{tikzpicture}
  \foreach \i in {2,...,6}{
  \path(\i,0) node[below]{\i};
  }
  \foreach \j in {0,...,3}{
  \path(2,\j) node[left]{\j};
  }
  \draw[step=1cm,color=gray] (2,0) grid (6,3);
  \draw[->,line width=2pt,font=\footnotesize] (3,3) -- node[left]{1} (3,2) -- node[above]{$1-q^2$} (4,2) -- node[above]{$\frac{1-q^2}{1-q^3}$} (5,2) -- node[above]{$\frac{1-q}{1-q^2}$} (6,2);
\end{tikzpicture}
\begin{tikzpicture}
  \foreach \i in {2,...,6}{
  \path(\i,0) node[below]{\i};
  }
  \foreach \j in {0,...,3}{
  \path(2,\j) node[left]{\j};
  }
  \draw[step=1cm,color=gray] (2,0) grid (6,3);
  \draw[->,line width=2pt,font=\footnotesize] (3,3) -- node[left]{1} (3,2) -- node[above]{$1-q^2$} (4,2) -- node[above]{$\frac{1-q^2}{1-q^3}$} (5,2) -- node[left]{$1-\frac{1-q}{1-q^2}$} (5,1) -- node[above]{1} (6,1);
\end{tikzpicture}
\begin{tikzpicture}
  \foreach \i in {2,...,6}{
  \path(\i,0) node[below]{\i};
  }
  \foreach \j in {0,...,3}{
  \path(2,\j) node[left]{\j};
  }
  \draw[step=1cm,color=gray] (2,0) grid (6,3);
  \draw[->,line width=2pt,font=\footnotesize] (3,3) -- node[left]{1} (3,2) -- node[above]{$1-q^2$} (4,2) -- node[left]{$1-\frac{1-q^2}{1-q^3}$} (4,1) -- node[above]{$1$} (5,1) -- node[above]{1} (6,1);
\end{tikzpicture}
\begin{tikzpicture}
  \foreach \i in {2,...,6}{
  \path(\i,0) node[below]{\i};
  }
  \foreach \j in {0,...,3}{
  \path(2,\j) node[left]{\j};
  }
  \draw[step=1cm,color=gray] (2,0) grid (6,3);
  \draw[->,line width=2pt,font=\footnotesize] (3,3) -- node[left]{1} (3,2) -- node[left]{$q^2$} (3,1) -- node[above]{$1$} (4,1) -- node[above]{$1$} (5,1) -- node[above]{1} (6,1);
\end{tikzpicture}
\end{center}
\caption{The four possible output tableaux $\tilde P$, their weights $I_3(P,\tilde P)$, and the corresponding
insertion paths, with edge weights indicated, for $k=3$ and $P$ given by \eqref{P}.}
\label{f:qlattice}
\end{figure}

The $q$-insertion algorithm can be applied to a word $w=w_1\ldots w_n\in[l]^n$,
starting with an empty tableau and successively inserting the numbers $w_1,w_2,\ldots,w_n$,
multiplying the weights along each possible sequence of output tableaux $P(1),\ldots,P(n)=P$
to obtain a distribution of weights $\phi_w(P,Q)$ on $\vart_l\times\vars_n$.  More precisely, we 
define $\phi_w(P,Q)$ recursively as follows.  Set
$$\phi_k(P,Q)=\begin{cases} 1 & \mbox{ if }P=k\mbox{ and }Q=1\\ 0 & \mbox{ otherwise.}\end{cases}$$
For $w\in [l]^n$ and $(\tilde P,\tilde Q)\in\vart_l\times\vars_{n+1}$ with $\sh \tilde P=\sh \tilde Q$, define
$$\phi_{wk}(\tilde P,\tilde Q)=\sum \phi_w(P,Q) I_k(P,\tilde P) ,$$
where the sum is over $(P,Q)\in \vart_l\times\vars_{n}$ with $\sh P=\sh Q$.

We conclude this section by giving a more algorithmic description of the $q$-insertion algorithm.
For this it is convenient to assume $0\le q<1$ and describe it using probabilistic language,
although it will be clear how to modify this using the language of `weights' in the general case.
For reference, we begin with an algorithmic description of the usual column insertion algorithm.
Denote the input word by $w\in[l]^n$.
\begin{enumerate}
  \item Set $i\leftarrow1$ and $(P,Q)=(\emptyset,\emptyset)$.
  \item Set $k\leftarrow w_i$ and $j\leftarrow k$.
  \item If ${\lambda^{k-1}_{j-1}=\lambda^k_j}$ and $j>1$ then set $j\leftarrow j-1$; otherwise $k$ displaces the first number $s$ in $j$th row of the tableau that is larger than $k$ ($s=\infty$ and $k$ is appended at the end of the row if no such number exists) and set $k\leftarrow s$.
  \item If $k=\infty$ then append $i$ to $Q$ such that $P$ and $Q$ have the same shape, set $i\leftarrow i+1$ and go to step (2); otherwise go to step (3).
\end{enumerate}
The $q$-insertion algorithim is defined as follows.
We adopt here the following convention: for $i>0$, let
\begin{align*}
  q^{\lambda^{i-1}_0-\lambda^i_1}=q^{\lambda^i_0-\lambda^i_1}=q^{\lambda^i_0-\lambda^{i-1}_0}=q^{\lambda^i_i-\lambda^i_{i+1}}=q^{\lambda^i_i-\lambda^{i-1}_i}=q^{\lambda^{i-1}_i-\lambda^i_{i+1}}=0.
\end{align*}
This convention is used for covering boundary conditions in general arguments.  It is only used in the following description of the 
$q$-insertion algorithm as well as in Section \ref{s:proofp}.  Otherwise the undefined $\lambda^i_j$ for $j>i$ or $j=0$ are 
taken to be zero.
\begin{enumerate}
  \item Set $i\leftarrow1$ and $(P,Q)=(\emptyset,\emptyset)$.
  \item Set $k\leftarrow w_i$, $j\leftarrow k$, $d\leftarrow0$ and $a_e(m,n)\leftarrow f_e(m,n)\ \forall e\in\{0,1\},1\le n\le m$.
  \item With probability $1-a_d(k,j)$ set $j\leftarrow j-1$ and $d\leftarrow0$; otherwise $k$ displaces the first number s in $j$th row of the tableau that is larger than $k$ ($s=\infty$ and append $k$ at the end of $j$th row if no such number exists) and set $k\leftarrow s$ and $d\leftarrow1$.
  \item If $k=\infty$ then append $i$ to $Q$ such that $P$ and $Q$ have the same shape, set $i\leftarrow i+1$ and go to step (2); otherwise go to step (3).
\end{enumerate}
As is obvious, when $q=0$ it reduces to the usual column insertion algorithm.

\section{Main result}

The weights $\phi_w(P,Q)$ are quite complicated.  The main result of this paper is that a remarkable
simplification occurs when we average over the set of words.  Before stating the result, we first introduce two
more functions on tableaux and explain their connection to $q$-Whittaker functions and Macdonald polynomials.
Denote the $q$-Pochhammer symbol by
$$(n)_q:=(q;q)_n=(1-q)\dots(1-q^n),$$ with the conventions $(n)_0=(0)_q=1$, 
and the $q$-binomial coefficients by
\begin{align*}
  \left[\begin{matrix}n\\m\end{matrix}\right]_q=\frac{(n)_q}{(m)_q(n-m)_q}.
\end{align*}
For $P\in\vart_l$ with $\sh P^i=\lambda^i$, $1\le i\le l$, writing $\lambda=\lambda^l$, define
\begin{align*}
  \kappa(P)&=\frac{\prod_{j=2}^{l-1}\prod_{i=1}^{j-1}(\lambda^j_i-\lambda^j_{i+1})_q}{\prod_{j=1}^{l-1}\prod_{i=1}^j(\lambda^j_i-\lambda^{j+1}_{i+1})_q(\lambda^{j+1}_i-\lambda^j_i)_q}\\
  &=\Delta_l(\lambda)^{-1}
  \prod_{1\le j<i\le l}\left[\begin{matrix}\lambda^i_j-\lambda^i_{j+1}\\\lambda^i_j-\lambda^{i-1}_j\end{matrix}\right]_q,
\end{align*}
where $$\Delta_l(\lambda)=\prod_{i=1}^{l-1}(\lambda_i-\lambda_{i+1})_q.$$
For $Q\in\vars_n$ with $\sh Q^i=\mu^i$, $1\le i\le n$, define
$$\rho(Q)=\prod_{1\le i\le j:\ \mu^i_j-\mu^{i-1}_j=1} (1-q^{\mu^i_j-\mu^i_{j+1}}).$$
The functions $\kappa$ and $\rho$ are simply related as follows.  Suppose that $l\ge n$ and $P$ has distinct entries
$i_1<i_2<\cdots<i_n$.  Denote by $\hat P\in\vars_n$ the standard tableau obtained by replacing the entry $i_k$
by $k$, for each $k=1,\ldots,n$.  Then
\begin{equation}\label{kapparho}
\kappa(P)=\frac{\rho(\hat P)}{(1-q)^n\Delta_l(\lambda)} .
\end{equation}
Indeed, using the simple identities,
\begin{align}\label{eq:simpbin}
  \left[\begin{matrix}a\\0\end{matrix}\right]_q=1,\qquad\qquad\left[\begin{matrix}a\\1\end{matrix}\right]_q=\frac{1-q^a}{1-q},
\end{align}
we have
  \begin{align*}
    \kappa( P)&=\Delta_l(\lambda)^{-1}\prod_{\substack{1\le i<j\le l\\\lambda^i_j-\lambda^{i-1}_j=1}}\left[\begin{matrix}\lambda^i_j-\lambda^i_{j+1}\\\lambda^i_j-\lambda^{i-1}_j\end{matrix}\right]_q\prod_{\substack{1\le i<j\le l\\\lambda^i_j-\lambda^{i-1}_j=0}}\left[\begin{matrix}\lambda^i_j-\lambda^i_{j+1}\\\lambda^i_j-\lambda^{i-1}_j\end{matrix}\right]_q\\
      &=\Delta_l(\lambda)^{-1}\prod_{\substack{1\le i<j\le l\\\lambda^i_j-\lambda^{i-1}_j=1}}\frac{1-q^{\lambda^i_j-\lambda^i_{j+1}}}{1-q}\prod_{\substack{1\le i=j\le l\\\lambda^i_j-\lambda^{i-1}_j=1}}\frac{1-q}{1-q}\\
    &=\frac{\rho(\hat P)}{(1-q)^n\Delta_l(\lambda)} .
    \end{align*}   

The functions $\kappa$ and $\rho$ are closely related to $q$-Whittaker functions~\cite{ruijsenaars,etingof,glo1,glo2}.  
Denote by  $\Omega^l$ the set of partitions with at most $l$ parts.  
The $q$-Whittaker function with parameter $a\in\C^l$
is a function on $\Omega^l$ defined by
\begin{align}\label{eq:Kiskernel}
  \Psi_a(\lambda)=\sum_{P\in\vart_l:\ \sh P=\lambda} a^P \kappa(P).
\end{align}
In~\cite{glo3} it is shown that these functions
are given in terms of the Macdonald polynomials $P_\lambda(x;q,t)$ as
\begin{equation}  \label{whi-mac}
\Psi_a(\lambda)= \Delta_l(\lambda)^{-1} P_\lambda(a;q,0).
\end{equation}

From this it follows that
\begin{equation}\label{eq:zetaq1}
\Psi_a(\lambda)= \Delta_l(\lambda)^{-1} \sum_\mu k_{\lambda\mu}(q) m_\mu(a)
\end{equation}
where $m_\mu$ denote the monomial symmetric functions and
\begin{equation}\label{k}
k_{\lambda\mu}(q)=\Delta_l(\lambda) \sum_{\sh P=\lambda, \ty P=\mu} \kappa(P)
= \sum_\nu K_{\lambda\nu}(q,0)K_{\nu\mu},
\end{equation}
where $K_{\lambda\mu}(q,t)$ are the two-variable Kostka polynomials~\cite{mac}.
We recall that $K_{\lambda\nu}(q,0)=K_{\lambda'\nu'}(0,q)=K_{\lambda'\nu'}(q)$, 
where $K_{\lambda\mu}(t)=K_{\lambda\mu}(0,t)$ are the single-variable Kostka polynomials.  
For an extensive survey of the various properties and interpretations of these polynomials, see~\cite{kirillov}.
When $q=0$, $\kappa(P)\equiv 1$ and $k_{\lambda\mu}(0)$ is equal to the Kostka number $K_{\lambda\mu}$, 
which is the number of tableaux with shape $\lambda$ and type $\mu$.  
In this case, $\Psi_a(\lambda)$ is given by the Schur polynomial 
$$\Psi_a(\lambda)=s_\lambda(a)=\sum_\mu K_{\lambda\mu} m_\mu(a).$$

We will also consider the following functions:
$$f^\lambda(q) = \sum_{Q\in\vars_n:\ \sh Q=\lambda}\rho(Q).$$
Note that $f^\lambda(0)=f^\lambda$, the number of standard tableaux with shape $\lambda$.
The relation between $f^\lambda(q)$ and the Whittaker functions $\Psi_a$ is given by the following proposition,
which is a straightforward consequence of \eqref{kapparho}.  Define
$$\Delta(\lambda)=\prod_{i=1}^{l(\lambda)}(\lambda_i-\lambda_{i+1})_q ,$$
where $l(\lambda)$ denotes the number of parts in $\lambda$.
\begin{prop}\label{limpsi} For each $\lambda\vdash n$,
$$\lim_{l\to\infty} \Psi_{(1/l)^l}(\lambda) = \frac{f^\lambda(q)}{n! (1-q)^n\Delta(\lambda)}.$$
\end{prop}
It follows, using
$$\lim_{l\to\infty} s_\lambda((1/l)^l) = f^\lambda / n!,$$
that $f^\lambda(q)$ is also given, for $\lambda\vdash n$, by
$$f^\lambda(q)=(1-q)^n \sum_\mu K_{\lambda\mu}(q,0) f^\mu .$$
To understand this in terms of specializations, recall that the exponential specialization $\mbox{ex}_1$ is the homomorphism 
defined on the ring of symmetric functions by $\mbox{ex}_1(p_n)=\delta_{n1}$, where $p_n$ are the elementary power sums 
(see, for example,~\cite[\S 7.8]{stanley}).  It follows from the above proposition (or can be seen directly) that 
$$f^\lambda(q) = n!(1-q)^n \mbox{ex}_1 \left( P_\lambda(q,0) \right).$$

The $q$-Whittaker functions $\Psi_a$ are eigenfunctions of Ruijsenaars' relativistic Toda difference 
operators~\cite{ruijsenaars, ruijsenaars1,etingof,glo1}.
In particular, 
\begin{equation}\label{ee}
L\Psi_a= \left(\sum_i a_i\right) \Psi_a,
\end{equation}
 where $L$ is the kernel operator defined by
$$L(\lambda,\mu)=\begin{cases} c_i(\lambda)&\mbox{ if $\mu=\lambda+e_i$ for some $1\le i\le l$}, \\
0 & \mbox{ otherwise,}\end{cases} $$
and
$$c_i(\lambda)=\begin{cases} 1-q^{\lambda_i-\lambda_{i+1}+1} & \mbox{ for } 1\le i<l,\\ 1 &  \mbox{ for } i=l.\end{cases}$$

Our main result is the following.
\begin{thm}\label{main}  
Let $(P,Q)\in \vart_l\times\vars_{n}$ with $\sh P=\sh Q=\lambda$.  Then
\begin{equation}\label{eq:main}
\sum_{w\in[l]^n} \phi_w(P,Q) = (\lambda_l)_q^{-1} \kappa(P) \rho(Q).
\end{equation}
\end{thm}
We note the following immediate extension of this identity which is useful for applications.
The type of a word $w$ is the composition $\mu=(\mu_1,\mu_2,\ldots)$ where $\mu_i$ is the number of $i$'s in $w$.
For $a=(a_1,\ldots,a_l)$ and $\mu$ a composition, write $a^\mu=a_1^{\mu_1}\ldots a_l^{\mu_l}$;
for $w\in[l]^n$ and $P\in\vart_l$, write $a^w=a^{\ty(w)}$ and $a^P=a^{\ty P}$.  Now, since $\phi_w(P,Q)=0$
unless $\ty P=\ty (w)$, we can write
\begin{equation}\label{main-a}
\sum_{w\in[l]^n} a^w \phi_w(P,Q) = (\lambda_l)_q^{-1} a^P \kappa(P) \rho(Q).
\end{equation}

Summing \eqref{main-a} over $P$ and $Q$ gives
$$\sum_{(P,Q)\in\vart_l\times\vars_n:\sh P=\sh Q=\lambda} \sum_{w\in[l]^n} a^w \phi_w(P,Q)
= (\lambda_l)_q^{-1} \Psi_a(\lambda) f^\lambda(q) .$$
Note that this implies the Cauchy-Littlewood type identity
$$\sum_{\lambda\vdash n} (\lambda_l)_q^{-1} \Psi_a(\lambda) f^\lambda(q) =\left(\sum_i a_i\right)^n.$$
Theorem~\ref{main} also yields some combinatorial formulas.
\begin{cly}  Let $\lambda,\mu\vdash n$ with at most $l$ parts, and let $Q$ be a standard tableau with shape $\lambda$.  
Then
$$P_\lambda(a;q,0) =  \sum_{w\in [l]^n} H_Q(w) a^w$$
and
$$k_{\lambda\mu}(q) =  \sum_{w\in [l]^n:\ \ty(w)=\mu} H_Q(w),$$
where
$$H_Q(w)= \dfrac{\Delta(\lambda)}{\rho(Q)} \sum_{P} \phi_w(P,Q).$$
Similarly, for any fixed $P\in\vart_l$ with shape $\lambda\vdash n$,
$$f^\lambda(q)=  \sum_{w\in [l]^n} G_P(w),$$
where
$$G_P(w)=\dfrac{(\lambda_l)_q}{\kappa(P)} \sum_{Q} \phi_w(P,Q).$$
Taking $P$ to be standard with shape $\lambda\vdash n$, this last formula becomes
$$\mbox{ex}_1 \left( P_\lambda(q,0) \right) = \frac1{n!} \sum_{\sigma\in S_n} H_P(\sigma),$$
where the sum is over permutations and $H_P(\sigma)$ indicates the function $H_P$ evaluated at the word
$\sigma^{-1}(1)\ldots \sigma^{-1}(n) $.
\end{cly}
When $q=0$, $H_S(w)$ equals 1 if the $Q$-tableau obtained by applying the Robinson-Schensted algorithm 
with column insertion to $w$ is $S$, and 0 otherwise; similarly, $G_T(w)$ equals 1 if the $P$-tableau obtained by applying 
the Robinson-Schensted algorithm with column insertion to $w$ is $T$, and 0 otherwise. 
The functions $G_T$ and $H_S$ thus generalise the notions of $P$-equivalence and $Q$-equivalence, 
or Knuth and dual Knuth equivalence, for the Robinson-Schensted algorithm with column insertion 
(see, for example, \cite[Chapter 2 and \S A.3]{fulton}).

The key ingredient in the proof of Theorem~\ref{main} is the following intertwining relation.  
Define kernel operators $K$ and $M$ by
$$K(\lambda,P)= a^P \kappa(P) \ind_{\sh P=\lambda},\qquad M(P,\tilde P)=\sum_{k=1}^l a_k I_k(P,\tilde P).$$
\begin{prop}\label{int} The following intertwining relation holds:
  \begin{align}\label{eq:int}
    KM=LK
  \end{align}
\end{prop}
We remark that \eqref{eq:int} immediately yields the eigenvalue equation \eqref{ee}.

\section{Stochastic evolutions}

If $0\le q<1$ and $a\in \real_+^l$ with $\sum_i a_i=1$, then 
\begin{equation}\label{pq}
\sum_{w\in[l]^n} a^w \phi_w(P,Q) = (\lambda_l)_q^{-1} a^P \kappa(P) \rho(Q)
\end{equation}
defines a probability measure on $\vart_l\times\vars_n$, which can be interpreted
as the distribution of the pair of tableaux obtained when one applies the randomised
insertion algotihm to a random word $w_1\ldots w_n$ with each $w_i$ chosen independently 
at random from $[l]$ according to the probabilities $a_1,\ldots,a_l$.  If we denote by $\mathcal L(m)$
the shape of the tableau obtained after inserting the first $m$ entries $w_1\ldots w_m$ then, 
given the interpretation of $Q$ as a recording tableau, we conclude by summing \eqref{pq}
over $P$ that the sequence of shapes $\mathcal L(1),\ldots,\mathcal L(n)$ is distributed 
according to
$$\prob(\mathcal L(1)=\mu^1,\ldots,\mathcal L(n)=\mu^n)=\Psi_a(\mu^n) \rho(Q),$$
where $Q\in\vars_n$ is defined by $\sh Q^i=\mu^i$, $i=1,\ldots,n$.  But this can be written as
$$\prob(\mathcal L(1)=\mu^1,\ldots,\mathcal L(n)=\mu^n)=
\prod_{i=1}^n \frac{\Psi_a(\mu^i)}{\Psi_a(\mu^{i-1})} L(\mu^{i-1},\mu^i).$$
Since $n$ is arbitrary, we immediately conclude the following.
Write $\mu\nearrow\lambda$ if $\lambda$ is obtained from $\mu$ by adding a single box.
\begin{thm}\label{mf}
When applying the randomised insertion algorithm to a random word $w_1w_2\ldots$ 
with each $w_i$ chosen independently at random from $[l]$ according to the probabilities $a_1,\ldots,a_l$
the sequence of tableaux $\mathcal P(n),\ n\ge 0$ obtained evolves as a Markov chain in $\vart_l$
with transition probabilities 
$$M(P,\tilde P)=\sum_{k=1}^n a_k I_k(P,\tilde P).$$
The sequence of shapes $\mathcal L(n)=\sh \mathcal P(n)$ evolves as a Markov chain in 
$\Omega^l$ with transition probabilities
$$p(\mu,\lambda)=\frac{\Psi_a(\lambda)}{\Psi_a(\mu)} L(\mu,\lambda) \ind_{\mu\nearrow\lambda}.$$
The conditional law of $\mathcal P(n)$, given
$\{ \mathcal L(1),\ldots,\mathcal L(n);\ \mathcal L(n)=\lambda \}$, is
$$\prob(\mathcal P(n)=P|\ \mathcal L(1),\ldots,\mathcal L(n);\ \mathcal L(n)=\lambda)=
\frac{K(\lambda,P)}{\Psi_a(\lambda)} .$$
The conditional law of $\ty \mathcal P(n)$, given
$\{ \mathcal L(1),\ldots,\mathcal L(n);\ \mathcal L(n)=\lambda \}$, is
$$\prob(\ty \mathcal P(n)=\mu|\ \mathcal L(1),\ldots,\mathcal L(n);\ \mathcal L(n)=\lambda)=
\frac{a^\mu k_{\lambda\mu}(q)}{\Psi_a(\lambda)} .$$
The distribution of $\mathcal L(n)$ is given by
$$\nu(\lambda):=\prob(\mathcal L(n)=\lambda)=(\lambda_l)_q^{-1} \Psi_a(\lambda) f^\lambda(q) .$$
\end{thm}
The probability distribution $\nu$ is a particular specialisation (and restriction to $\lambda\vdash n$)
of the Macdonald measures introduced by Forrester and Rains~\cite{fr}, see also~\cite{bc}.
When $q=0$, the above theorem reduces to the fact~\cite{noc1} that, when applying the usual column insertion
algorithm to a random word with probabilities $a_1,\ldots,a_l$, the shape of the tableau evolves
as a Markov chain with transition probabilities
$$p(\mu,\lambda)=\frac{s_\lambda(a)}{s_\mu(a)} \ind_{\mu\nearrow\lambda}.$$
If $a_1>a_2>\cdots>a_l$ this Markov chain can be interpreted as a random walk in $\N^l$ with transition probabilities
$$r(\mu,\lambda)=a^{\lambda-\mu}  \ind_{\mu\nearrow\lambda}$$
conditioned never to exit the Weyl chamber $\{\lambda\in\N^l:\ \lambda_1\ge\cdots\ge\lambda_l\}$,
which can be identified with $\Omega^l$.  This result, which relates to the representation theory of $\mathfrak{gl}_l$, 
has been generalised to arbitrary complex semisimple Lie algebras in~\cite{bbo,llp}.  
For earlier related work on the asymptotics of longest monotone subsequences in random words, see~\cite{tw1}.
When $q\to1$ the $q$-Whittaker functions converge with appropriate rescaling to $\gl_l$-Whittaker functions \cite{glo2}, 
and the above theorem should re-scale to the main result of the paper~\cite{noc12}, which relates a continuous-time
version of the geometric RSK correspondence introduced by A.N. Kirillov~\cite{kirillov1}, with Brownian motion 
as input, to the open quantum Toda chain with $l$ particles.  In this scaling limit, the $q$-insertion 
algorithm should converge in an appropriate sense to the continuous-time version of the 
geometric RSK mapping considered in~\cite{noc12}, which is deterministic.  
The results of~\cite{noc12} have been generalised in~\cite{ch} (see also \cite{bbo1}) to arbitrary complex semisimple
Lie algebras.  It is natural to expect the results of the present paper to admit a similar generalisation.

\begin{exa}
  The rank-1 case ($l=2$) of Theorem~\ref{mf} is discussed in \cite{noc12b}. 
  Setting $\mathcal L^i(n)=\sh\mathcal P^i(n)$, the evolution on tableaux in this case is driven 
  by the process $Y(n)=\mathcal L^1_1(n)-\mathcal L^2_1(n)$, $n\ge 0$,
  which (setting $p=a_1$) is a birth and death process as illustrated in Figure \ref{fig:bdy}.
\begin{figure}[h!]
  \begin{center}
    \begin{tikzpicture}[font=\footnotesize,scale=2,>=stealth]
      \tikzstyle{every node}=[draw,shape=circle,minimum size=35pt];
      \path(0,0) node (v0) {$0$};
      \path(1,0) node (v1) {$1$};
      \path(2,0) node (v2) {$2$};
      \path(3,0) node (vy0) {$y-1$};
      \path(4,0) node (vy1) {$y$};
      \path(5,0) node (vy2) {$y+1$};
      \path(5.7,0) coordinate (end);

      \tikzstyle{every node}=[];
      \draw [->,thick] (v0) .. controls +(-135:.7cm) and +(135:.7cm) .. node[left]{$1-p$} (v0);

      \draw [->,thick] (v0) .. controls +(45:.5cm) and +(135:.5cm) .. node[above]{$p$} (v1);
      \draw [->,thick] (v1) .. controls +(45:.5cm) and +(135:.5cm) .. node[above]{$p$} (v2);
      \draw [color=white] (v2) -- node[color=black]{$\dots$} (vy0);
      \draw [->,thick] (vy0) .. controls +(45:.5cm) and +(135:.5cm) .. node[above]{$p$} (vy1);
      \draw [->,thick] (vy1) .. controls +(45:.5cm) and +(135:.5cm) .. node[above]{$p$} (vy2);
      \draw [color=white] (vy2) -- node[color=black]{$\dots$} (end);

      \draw [->,thick] (v1) .. controls +(-135:.5cm) and +(-45:.5cm) .. node[below]{$(1-p)(1-q)$} (v0);
      \draw [->,thick] (v2) .. controls +(-135:.5cm) and +(-45:.5cm) .. node[below]{$(1-p)(1-q^2)$} (v1);
      \draw [->,thick] (vy1) .. controls +(-135:.5cm) and +(-45:.5cm) .. node[below]{$(1-p)(1-q^y)$} (vy0);
      \draw [->,thick] (vy2) .. controls +(-135:.5cm) and +(-45:.5cm) .. node[below]{$(1-p)(1-q^{y+1})$} (vy1);

      \draw [->,thick] (v1) .. controls +(125:.7cm) and +(55:.7cm) .. node[above]{$(1-p)q$} (v1);
      \draw [->,thick] (v2) .. controls +(125:.7cm) and +(55:.7cm) .. node[above]{$(1-p)q^2$} (v2);
      \draw [->,thick] (vy0) .. controls +(125:.7cm) and +(55:.7cm) .. node[above]{$(1-p)q^{y-1}$} (vy0);
      \draw [->,thick] (vy1) .. controls +(125:.7cm) and +(55:.7cm) .. node[above]{$(1-p)q^{y}$} (vy1);
      \draw [->,thick] (vy2) .. controls +(125:.7cm) and +(55:.7cm) .. node[above]{$(1-p)q^{y+1}$} (vy2);
    \end{tikzpicture}
  \end{center}
  \caption{The birth-and-death process $Y$}
  \label{fig:bdy}
\end{figure}
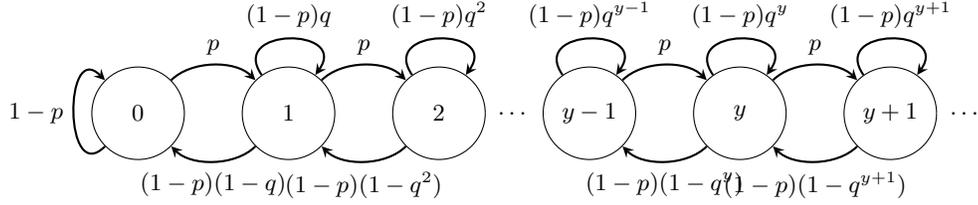
\end{exa}
\begin{exa}\label{e:rank2}
  When $l=3$ the algorithm is more complicated than in the
  $l=2$ case because the push-or-bump probability $f_1(3,2)$ appears. 
  In this case the algorithm with random input is described as follows (cf. Example 2).
  In the following, w.p. means ``with probability''.
  \begin{itemize}
    \item w.p. $a_1$, insert 1 to row 1, pushing 2's and 3's in row 1
    \item w.p. $a_2$, insert 2
      \begin{itemize}
	\item w.p. $1-q^{\lambda^1_1-\lambda^2_2}$, the 2 is inserted to row 2 and the displaced 3 is either pushed or bumped
	\begin{itemize}
	  \item w.p. $(1-q^{\lambda^2_1-\lambda^3_2})/(1-q^{\lambda^2_1-\lambda^2_2})$ the displaced 3 is pushed in row 2
	  \item w.p. $1-(1-q^{\lambda^2_1-\lambda^3_2})/(1-q^{\lambda^2_1-\lambda^2_2})$ the displaced 3 is bumped to row 1
	\end{itemize}
	\item w.p. $q^{\lambda^1_1-\lambda^2_2}$, the 2 is inserted to row 1 and it pushes 3's in row 1
      \end{itemize}
    \item w.p. $a_3$, insert 3
      \begin{itemize}
	\item w.p. $1-q^{\lambda^2_2-\lambda^3_3}$, the 3 is inserted to row 3
	\item w.p. $q^{\lambda^2_2-\lambda^3_3}(1-q^{\lambda^2_1-\lambda^3_2})$ the 3 is inserted to row 2
	\item w.p. $q^{\lambda^2_2-\lambda^3_3}q^{\lambda^2_1-\lambda^3_2}$ the 3 is inserted to row 1
      \end{itemize}
  \end{itemize}
\end{exa}

The $q$-insertion algorithm applied to a random word is closely related to the $q$-TASEP
interacting particle system.  This is a variation of the totally asymmetric simple exclusion process 
(TASEP) which was introduced (in the language of $q$-bosons) and shown to be integrable
by Sasamoto and Wadati \cite{sw}, and recently related to $q$-Whittaker functions by Borodin and 
Corwin \cite{bc}.  The process is defined as follows.
There are $l$ particles on the integer lattice, and we denote their positions by $x_1>x_2>\dots>x_l$. 
Let $a_1,a_2,\dots,a_l\in\real_+$.  Without loss of generality we can assume $\sum_i a_i=1$.
The particles jump independently to the right by 1 with respective rates
\begin{align*}
  r_i=
  \begin{cases}
    a_1,&\text{ if }i=1;\\
    a_i(1-q^{x_{i-1}-x_i-1}),&\text{ otherwise.}
  \end{cases}
\end{align*}
Note that when $x_i+1=x_{i-1}$ the rate $r_i$ vanishes, thus enforcing the exclusion rule.
Now consider the tableau-valued Markov chain $\mathcal P(n)$, $n\ge 0$, defined as above
by applying the randomised insertion algorithm applied to a random word with probabilities $a_1,\ldots,a_l$.
Setting $\mathcal L^i(n)=\sh\mathcal P^i(n)$, we see that the process $X_1(n),\ldots,X_l(n),\ n\ge 0$
defined by $X_i(n)=\mathcal L^i_i(n)-i+1$ evolves as a Markov chain with state space 
$\{x\in\Inte^l:\ x_1>x_2>\cdots>x_l\}$ and transition
probabilities 
$$\pi(x,x+e_i)=r_i,\ i=1,\ldots,l\qquad\pi(x,x)=1-\sum_i r_i,$$ 
where $r_i$ are defined as above.  
In other words, it is a de-Poissonisation of the $q$-TASEP process.  
Denote the $q$-TASEP process by $\tilde X(t),\ t\ge 0$, started with step initial 
condition $\tilde X_i(0)=1-i$, $i=1,\ldots,l$; by Theorem~\ref{mf}, the law of the position of the last particle at 
time $t$ is given by
\begin{align}\label{qt}
  \prob(\tilde X_l(t)=m-l+1) &=\sum_{k\ge0}e^{-t}\frac{t^k}{k!}\sum_{\lambda\vdash k,\lambda^l=m}(\lambda_l)_q^{-1}
  \Psi_a(\lambda)f^\lambda(q)\nonumber \\
  &=e^{-t}\sum_{\lambda\in\Omega^l,\lambda_l=m}\frac{t^{|\lambda|}}{|\lambda|!}(\lambda_l)_q^{-1}\Psi_a(\lambda)f^\lambda(q).
\end{align}
In \cite{bc}, a continuous-time Markov chain on the set of tableaux $\vart_l$ (actually discrete Gelfand-Tsetlin patterns,
but this is equivalent)
was introduced.  It has the same fixed time marginals as the Poissonisation of the process $\mathcal P(n)$, 
although the dynamics are quite different.  It is also coupled in exactly the same way to the $q$-TASEP process 
and in the paper \cite{bc} an equivalent expression to \eqref{qt} is obtained via this coupling for the law of $\tilde X_l(t)$.
See also~\cite{bcs} for related recent work.

\section{Permutations}

If $l=n$ and $P\in\vars_n$ with $\sh P=\lambda$, then \eqref{kapparho} becomes
$$\kappa(P)=\frac{\rho(P)}{(1-q)^n\Delta_n(\lambda)} .$$  Using this, and the fact
that $\phi_w(P,Q)=0$ unless $\ty P=\ty w$, we immediately deduce from Theorem~\ref{main}
the following corollary.
\begin{cly}\label{perms} For $P,Q\in\vars_n$ with $\sh P=\sh Q=\lambda$, we have
  \begin{align}
\zeta_{P,Q}(q) := \sum_{\sigma\in S_n} \phi_\sigma(P,Q) = 
\frac{\rho(P)\rho(Q)}{(1-q)^n\Delta(\lambda)} .
    \label{eq:zetaq2}
  \end{align}
  \end{cly}
Summing over $P$ and $Q$ gives
$$\theta_\lambda(q) := \sum_{P,Q\in\vars_n:\sh P=\sh Q=\lambda} \zeta_{P,Q}(q) 
=\frac{f^\lambda(q)^2}{(1-q)^n\Delta(\lambda)}.$$
We note that $\sum_{\lambda\vdash n}\theta_\lambda(q)=n!$.  When $0\le q<1$, the probability 
measure on integer partitions defined by $\mu_q(\lambda)=\theta_\lambda(q)/n!$ gives the law of 
the shape of the tableaux obtained when one applied the randomised insertion algorithm to a
random permutation.  It would be interesting to understand the analogue in this setting
of the longest increasing subsequence problem~\cite{ad,bdj,okounkov}.

For any standard tableau $P$ with entries in $[n]$ and shape $\lambda$. Its weight $\rho(P)$ is a product of $n$ polynomials of the form of $(1-q^k)$ and hence $\rho(P)$ is divisible by $(1-q)^n$. On the other hand, considering the $i$th and $i+1$th row in $P$, each time $j$ a box is added in $i$th row, a factor $(1-q^d)$ - where $d$ is the difference between length of the corresponding
two rows at time $j$ - appears in $\rho(P)$.  For this difference $d$ to reach the value of $\lambda_i-\lambda_{i+1}$ eventually
(which it evidently does) all the factors $(1-q),(1-q^2),\dots,(1-q^{\lambda_i-\lambda_{i+1}})$ must appear at least once. 
It follows that $\rho(P)$ is also divisible by $\Delta(\lambda)$. 
Thus, $\zeta_{P, Q}(q)\in\Inte[q]$ for each pair $(P,Q)$ and $\theta_\lambda(q)\in\Inte[q]$ for each $\lambda$.

For any permutation $\sigma\in S_n$, denote by $(P(\sigma),Q(\sigma))$ the pair of tableaux after column inserting $\sigma$,
and set $F_\sigma(q)=\zeta_{P(\sigma),Q(\sigma)}(q)$.  When $n=2$, the polynomials $F_\sigma(q)$ and $\theta_\lambda(q)$
are given by
  \begin{align*}
    F_{12}(q)=1-q;&&F_{21}(q)=1+q.
  \end{align*}
  \begin{align*}
    \theta_{2}(q)=1+q;&&\theta_{1^2}(q)=1-q.
  \end{align*}
  When $n=3$, we have
  \begin{align*}
    F_{123}(q)=(1-q)^2;&&F_{132}(q)=1-q;&&F_{213}(q)=(1+q)(1-q^2);\\
    F_{231}(q)=1-q^2;&&F_{312}(q)=1-q^2;&&F_{321}(q)=(1+q)(1+q+q^2).
  \end{align*}
  \begin{align*}
    \theta_{3}(q)=(1+q)(1+q+q^2);&& \theta_{21}(q)=(1-q)(2+q)^2;&& \theta_{1^3}(q)=(1-q)^2.
  \end{align*}
  
The polynomials $F_\sigma(q)$ give an alternative interpretation of the probability measure
 $\mu_q$ as the distribution of the shape of the tableaux obtained when one applies the Robinson-Schensted 
column insertion algorithm to a permutation chosen at random according to the distribution $F_\sigma(q)/n!$.  

\section{Proofs}\label{s:proofs}
\subsection{Proof of Proposition \ref{int}}\label{s:proofp}
To prove \eqref{eq:int}, we take advantage of the recursive structure of the $q$-Whittaker functions. 
Define $\hat\kappa$ on $\Omega^l\times\Omega^{l-1}$ by
\begin{align*}
  \hat\kappa(\lambda^l,\lambda^{l-1})=\frac{\prod_{i=1}^{l-2}(\lambda^{l-1}_i-\lambda^{l-1}_{i+1})_q}{\prod_{i=1}^{l-1}(\lambda^{l-1}_i-\lambda^l_{i+1})_q(\lambda^l_i-\lambda^{l-1}_i)_q}
\end{align*}
and set 
$$T=\{(\lambda^l,\lambda^{l-1})\in\Omega^l\times\Omega^{l-1}:\lambda^{l-1}\prec\lambda^l\},$$
where we write $\lambda\prec\mu$ if $\mu_{i+1}\le \lambda_i\le \mu_i$ for each $i$.

We begin by verifying the simpler intertwining relation:
  \begin{align}\label{eq:recint}
    \hat K\hat M=L\hat K,
  \end{align}
  where $\hat M:T\times T\to\real_{\ge0}$ and $\hat K:T\to\real_{\ge0}$ are defined as follows.
  \begin{align*}
    \hat M((\lambda^l,\lambda^{l-1}),&(\lambda^l+\ve_k,\lambda^{l-1}))=a_l(1-q^{\lambda^{l-1}_{k-1}-\lambda^l_k})\prod_{i=k}^{l-1}q^{\lambda^{l-1}_i-\lambda^l_{i+1}},\quad 1\le k\le l;
  \end{align*}
  \begin{align*}
    \hat M((\lambda^l,\lambda^{l-1}),(\lambda^l+\ve_k,\lambda^{l-1}+\ve_k))
    =\frac{(1-q^{\lambda^{l-1}_k-\lambda^{l-1}_{k+1}+1})(1-q^{\lambda^{l-1}_{k-1}-\lambda^l_k})}{1-q^{\lambda^{l-1}_{k-1}-\lambda^{l-1}_k}},\\\quad 1\le k\le l-1;
  \end{align*}
  \begin{align*}
    \hat M((\lambda^l,&\lambda^{l-1}),(\lambda^l+\ve_k,\lambda^{l-1}+\ve_m))\\
    &=\frac{(1-q^{\lambda^{l-1}_m-\lambda^{l-1}_{m+1}+1})(1-q^{\lambda^l_m-\lambda^{l-1}_m})}{1-q^{\lambda^{l-1}_{m-1}-\lambda^{l-1}_m}}(1-q^{\lambda^{l-1}_{k-1}-\lambda^l_k})\prod_{i=k+1}^mq^{\lambda^{l-1}_{i-1}-\lambda^l_i},\\
    &\qquad\qquad\qquad\qquad\qquad\qquad\qquad\qquad\qquad\qquad1\le k<m\le l-1.
  \end{align*}
  \begin{align*}
  \hat K(\lambda^l,(\tilde\lambda^l,\lambda^{l-1}))
  &=a^{\sum_{i=1}^l\lambda^l_i-\sum_{i=1}^{l-1}\lambda^{l-1}_i}\hat\kappa(\lambda^l,\lambda^{l-1}) 
  \ind_{\lambda^l=\tilde\lambda^l}.
  \end{align*}
With a slight abuse of notation we will write $\hat K(\lambda^l,\lambda^{l-1})$ as shorthand for
  $\hat K(\lambda^l,(\tilde\lambda^l,\lambda^{l-1}))$ since the support of latter is in $\{\lambda^l=\tilde\lambda^l\}$. 
  We'll do the same for kernel $K$.
 
 We will verify the recursive intertwining relation \eqref{eq:recint} directly.  The left hand side is given by
  \begin{align*}
    \hat K\hat M(\lambda^l,&(\lambda^l+\ve_k,\lambda^{l-1}))=\hat K(\lambda^l,\lambda^{l-1})\hat M( (\lambda^l,\lambda^{l-1}),(\lambda^l+\ve_k,\lambda^{l-1}))\\
    &+\hat K(\lambda^l,\lambda^{l-1}-\ve_k)\hat M( (\lambda^l,\lambda^{l-1}-\ve_k),(\lambda^l+\ve_k,\lambda^{l-1}))\mathbb I_{k\le l-1}\\
    &+\sum_{m=k+1}^{l-1}\hat K(\lambda^l,\lambda^{l-1}-\ve_m)\hat M( (\lambda^l,\lambda^{l-1}-\ve_m),(\lambda^l+\ve_k,\lambda^{l-1}))\mathbb I_{k\le l-2}.
  \end{align*}
  We calculate each term separately. Set $K'=a_l\hat K(\lambda^l,\lambda^{l-1})$.
  \begin{align*}
    \hat K(\lambda^l,\lambda^{l-1})\hat M( (\lambda^l,\lambda^{l-1}),(\lambda^l+\ve_k,\lambda^{l-1}))=K'(1-q^{\lambda^{l-1}_{k-1}-\lambda^l_k})\prod_{i=k}^{l-1}q^{\lambda^{l-1}_i-\lambda^l_{i+1}}.
  \end{align*}
  \begin{align*}
    \hat K(\lambda^l,&\lambda^{l-1}-\ve_k)\hat M( (\lambda^l,\lambda^{l-1}-\ve_k),(\lambda^l+\ve_k,\lambda^{l-1}))\\
    &=K'\frac{(1-q^{\lambda^{l-1}_{k-1}-\lambda^{l-1}_k+1})(1-q^{\lambda^{l-1}_k-\lambda^l_{k+1}})}{(1-q^{\lambda^{l-1}_k-\lambda^{l-1}_{k+1}})(1-q^{\lambda^l_k-\lambda^{l-1}_k+1})}\frac{(1-q^{\lambda^{l-1}_k-\lambda^{l-1}_{k+1}})(1-q^{\lambda^{l-1}_{k-1}-\lambda^l_k})}{1-q^{\lambda^{l-1}_{k-1}-\lambda^{l-1}_k+1}}\\
    &=K'(1-q^{\lambda^{l-1}_{k-1}-\lambda^l_k})\frac{1-q^{\lambda^{l-1}_k-\lambda^l_{k+1}}}{1-q^{\lambda^l_k-\lambda^{l-1}_k+1}}.
  \end{align*}
  \begin{align*}
    \sum_{m=k+1}^{l-1}\hat K(\lambda^l,&\lambda^{l-1}-\ve_m)\hat M( (\lambda^l,\lambda^{l-1}-\ve_m),(\lambda^l+\ve_k,\lambda^{l-1}))\\
    &=K'\sum_{m=k+1}^{l-1}\frac{(1-q^{\lambda^{l-1}_{m-1}-\lambda^{l-1}_m+1})(1-q^{\lambda^{l-1}_m-\lambda^l_{m+1}})}{(1-q^{\lambda^{l-1}_m-\lambda^{l-1}_{m+1}})(1-q^{\lambda^l_m-\lambda^{l-1}_m+1})}\\
    &\times\frac{(1-q^{\lambda^{l-1}_m-\lambda^{l-1}_{m+1}})(1-q^{\lambda^l_m-\lambda^{l-1}_m+1})}{1-q^{\lambda^{l-1}_{m-1}-\lambda^{l-1}_m+1}}(1-q^{\lambda^{l-1}_{k-1}-\lambda^l_k})\prod_{i=k+1}^mq^{\lambda^{l-1}_{i-1}-\lambda^l_i}\\
    &=K'(1-q^{\lambda^{l-1}_{k-1}-\lambda^l_k})\sum_{m=k+1}^{l-1}(1-q^{\lambda^{l-1}_m-\lambda^l_{m+1}})\prod_{i=k+1}^mq^{\lambda^{l-1}_{i-1}-\lambda^l_i}.
  \end{align*}
 The left hand side of \eqref{eq:recint} is thus given by
  \begin{align*}
    LHS&=K'(1-q^{\lambda^{l-1}_{k-1}-\lambda^l_k})\left(\prod_{i=k}^{l-1}q^{\lambda^{l-1}_i-\lambda^l_{i+1}}\right.\\
    &\left.+\sum_{m=k+1}^{l-1}(1-q^{\lambda^{l-1}_m-\lambda^l_{m+1}})\prod_{i=k+1}^mq^{\lambda^{l-1}_{i-1}-\lambda^l_i}\mathbb I_{k\le l-2}+\frac{1-q^{\lambda^{l-1}_k-\lambda^l_{k+1}}}{1-q^{\lambda^l_k-\lambda^{l-1}_k+1}}\mathbb I_{k\le l-1}\right)\\
    &\qquad\qquad\qquad\qquad=K'(1-q^{\lambda^{l-1}_{k-1}-\lambda^l_k})\frac{1-q^{\lambda^l_k-\lambda^l_{k+1}+1}}{1-q^{\lambda^l_k-\lambda^{l-1}_k+1}}.
  \end{align*}
  The right hand side is much easier to calculate:
  \begin{align*}
    L\hat K(\lambda^l,(\lambda^l+\ve_k,\lambda^{l-1}))=L(\lambda^l,\lambda^l+\ve_k)\hat K(\lambda^l+\ve_k,\lambda^{l-1})\\
    =K'(1-q^{\lambda^l_k-\lambda^l_{k+1}+1})\frac{1-q^{\lambda^{l-1}_{k-1}-\lambda^l_k}}{1-q^{\lambda^l_k-\lambda^{l-1}_k+1}},
  \end{align*}
as required.
  
We will now prove \eqref{eq:int} by induction on $l$.  When $l=2$, since $\hat M^2$ is the kernel for the whole tableau, 
the recursive intertwining relation \eqref{eq:recint} is equivalent to the full intertwining relation \eqref{eq:int}. 
Suppose the statement of the proposition holds for the rank-$(l-2)$ case, that is, for $l-1$.
From the definition of $K$ and $\hat K$ we have
    \begin{align}
      K^l(\lambda^l,\lambda^{1:l-1})=K^{l-1}(\lambda^{l-1},\lambda^{1:l-2})\hat K^l(\lambda^l,\lambda^{l-1}).
      \label{eq:krec}
    \end{align}
 By the recursive nature of definition of $\phi_w$, $M^l$ can be expressed in terms of $\hat M^l$, $M^{l-1}$ and $L^{l-1}$:
    \begin{align*}
      M^l(\lambda^{1:l},\tilde\lambda^{1:l})&=\ind_{\lambda^{l-1}=\tilde\lambda^{l-1}}\hat M^l( (\lambda^l,\lambda^{l-1}),(\tilde\lambda^l,\tilde\lambda^{l-1}))\\
      &+\ind_{\lambda^{l-1}\nearrow\tilde\lambda^{l-1}}\frac{M^{l-1}(\lambda^{1:l-1},\tilde\lambda^{1:l-1})}{L^{l-1}(\lambda^{l-1},\tilde\lambda^{l-1})}\hat M^l( (\lambda^l,\lambda^{l-1}),(\tilde\lambda^l,\tilde\lambda^{l-1})).
    \end{align*}
For partitions $\lambda,\mu$ write $\lambda\leadsto\mu$ to mean that either $\lambda=\mu$ or $\lambda\nearrow\mu$. 
Then
    \begin{align*}
      K^lM^l&(\lambda^l,(\tilde\lambda^l,\lambda^{1:l-1}))=\sum_{\tilde\lambda^{1:l-1}:\tilde\lambda^{l-1}\leadsto\lambda^{l-1}}K^l(\lambda^l,\tilde\lambda^{1:l-1})M^l( (\lambda^l,\tilde\lambda^{1:l-1}),(\tilde\lambda^l,\lambda^{1:l-1}))\\
      &=\sum_{\tilde\lambda^{1:l-1}}\hat K^l(\lambda^l,\tilde\lambda^{l-1})K^{l-1}(\tilde\lambda^{l-1},\tilde\lambda^{1:l-2})\Big(\ind_{\tilde\lambda^{l-1}=\lambda^{l-1}}\hat M^l( (\lambda^l,\tilde\lambda^{l-1}),(\tilde\lambda^l,\lambda^{l-1}))\\
      &+\ind_{\tilde\lambda^{l-1}\nearrow\lambda^{l-1}}\frac{M^{l-1}(\tilde\lambda^{1:l-1},\lambda^{1:l-1})}{L^{l-1}(\tilde\lambda^{l-1},\lambda^{l-1})}\hat M^l( (\lambda^l,\tilde\lambda^{l-1}),(\tilde\lambda^l,\lambda^{l-1}))\Big)\\
      &\qquad\qquad\qquad\qquad\qquad=:\ind_{\tilde\lambda^{l-1}=\lambda^{l-1}}\text I+\ind_{\tilde\lambda^{l-1}\nearrow\lambda^{l-1}}\text{II}.
    \end{align*}
    \begin{align*}
      \text{II}&=\sum_{\tilde\lambda^{l-1}\nearrow\lambda^{l-1}}\Bigg(\hat K^l(\lambda^l,\tilde\lambda^{l-1})\hat M^l( (\lambda^l,\tilde\lambda^{l-1}),(\tilde\lambda^l,\lambda^{l-1}))\\
      &\qquad\qquad\times\sum_{\tilde\lambda^{1:l-2}:\tilde\lambda^{l-2}\leadsto\lambda^{l-2}}K^{l-1}(\tilde\lambda^{l-1},\tilde\lambda^{1:l-2})\frac{M^{l-1}(\tilde\lambda^{1:l-1},\lambda^{1:l-1})}{L^{l-1}(\tilde\lambda^{l-1},\lambda^{l-1})}\Bigg)\\
      &=\sum_{\tilde\lambda^{l-1}\nearrow\lambda^{l-1}}\Bigg(\hat K^l(\lambda^l,\tilde\lambda^{l-1})\hat M^l( (\lambda^l,\tilde\lambda^{l-1}),(\tilde\lambda^l,\lambda^{l-1}))\\
      &\qquad\qquad\qquad\qquad\qquad\qquad\qquad\times\frac{K^{l-1}M^{l-1}(\tilde\lambda^{l-1},(\lambda^{l-1},\lambda^{1:l-2}))}{L^{l-1}(\tilde\lambda^{l-1},\lambda^{l-1})}\Bigg)\\
      &\xlongequal[\text{assumption}]{\text{induction}}\sum_{\tilde\lambda^{l-1}\nearrow\lambda^{l-1}}\Bigg(\hat K^l(\lambda^l,\tilde\lambda^{l-1})\hat M^l( (\lambda^l,\tilde\lambda^{l-1}),(\tilde\lambda^l,\lambda^{l-1}))\\
      &\qquad\qquad\qquad\qquad\qquad\qquad\qquad\times\frac{L^{l-1}K^{l-1}(\tilde\lambda^{l-1},(\lambda^{l-1},\lambda^{1:l-2}))}{L^{l-1}(\tilde\lambda^{l-1},\lambda^{l-1})}\Bigg)\\
      &=\sum_{\tilde\lambda^{l-1}\nearrow\lambda^{l-1}}\hat K^l(\lambda^l,\tilde\lambda^{l-1})\hat M^l( (\lambda^l,\tilde\lambda^{l-1}),(\tilde\lambda^l,\lambda^{l-1}))K^{l-1}(\lambda^{l-1},\lambda^{1:l-2})
    \end{align*}
    Due to the indicator, when $\tilde\lambda^{l-1}=\lambda^{l-1}$,
    \begin{align*}
      \text I=\hat K^l(\lambda^l,\tilde\lambda^{l-1})K^{l-1}(\lambda^{l-1},\lambda^{1:l-2})\hat M^l( (\lambda^l,\tilde\lambda^{l-1}),(\tilde\lambda^l,\lambda^{l-1})).
    \end{align*}
  Therefore 
  \begin{align*}
    K^lM^l&(\lambda^l,(\tilde\lambda^l,\lambda^{1:l-1}))\\
    &=\sum_{\tilde\lambda^{l-1}:\tilde\lambda^{l-1}\leadsto\lambda^{l-1}}\hat K^l(\lambda^l,\tilde\lambda^{l-1})K^{l-1}(\lambda^{l-1},\lambda^{1:l-2})\hat M^l( (\lambda^l,\tilde\lambda^{l-1}),(\tilde\lambda^l,\lambda^{l-1}))\\
    &=K^{l-1}(\lambda^{l-1},\lambda^{1:l-2})\hat K^l\hat M^l(\lambda^l,(\tilde\lambda^l,\lambda^{l-1}))\\
    &\overset{\eqref{eq:recint}}{=}K^{l-1}(\lambda^{l-1},\lambda^{1:l-2})L^l(\lambda^l,\tilde\lambda^l)\hat K^l(\tilde\lambda^l,\lambda^{l-1})\\
    &\overset{\eqref{eq:krec}}{=}L^l(\lambda^l,\tilde\lambda^l)K^l(\tilde\lambda^l,\lambda^{1:l-1})=LK(\lambda^l,(\tilde\lambda^l,\lambda^{1:l-1})),
  \end{align*}
  as required.

\subsection{Proof of Theorem \ref{main}}

We will prove the identity \eqref{main-a}, from which the statement of the theorem follows.
From the definition of $\phi_w$, for $(P,Q)\in\vart_l\times\vars_n$ such that $\sh P=\sh Q=\lambda$ and $\mu^i=\sh Q^i$
for $i=1,\ldots,n$, the left hand side of \eqref{main-a} can be written as
  \begin{align*}
    &\sum_{w\in[l]^n}a^w\phi_w(P,Q)\\
    &=\sum_{w\in[l]^n}\sum_{(P(i))_{i=1}^{n-1},\sh P(i)=\mu^i}a^wI_{w_1}(\emptyset,P(1))\dots I_{w_n}(P(n-1),P)\\
    &=\sum_{w\in[l]^n}\sum_{(P(i))_{i=1}^{n-1},\sh P(i)=\mu^i}(a_{w_1}I_{w_1}(\emptyset,P(1)))\dots(a_{w_n}I_{w_n}(P(n-1),P))\\
    &=\sum_{(P(i))_{i=1}^{n-1},\sh P(i)=\mu^i}\left(\sum_{w_1\in[l]}a_{w_1}I_{w_1}(\emptyset,P(1))\right)\dots\left(\sum_{w_n\in[l]}a_{w_n}I_{w_n}(P(n-1),P)\right)\\
    &=\sum_{(P(i))_{i=1}^{n-1},\sh P(i)=\mu^i}M(\emptyset,P(1))\dots M(P(n-1),P).
  \end{align*}
  
  On the right hand side, from the definition of $\rho(Q)$ and the intertwining relation \eqref{eq:int}, 
  \begin{align*}
    a^P&\kappa(P)\frac{\rho(Q)}{(\lambda_l)_q}=L(\emptyset,\mu^1)\dots L(\mu^{n-1},\lambda)K(\lambda,P)\\
    &=L(\emptyset,\mu^1)\dots L(\mu^{n-2},\mu^{n-1})LK(\mu^{n-1},P)\\
    &=L(\emptyset,\mu^1)\dots L(\mu^{n-2},\mu^{n-1})KM(\mu^{n-1},P)\\ 
    &=\sum_{P(n-1):\ \sh P(n-1)=\mu^{n-1}}\Bigg(L(\emptyset,\mu^1)\dots L(\mu^{n-2},\mu^{n-1})\\
    &\qquad\qquad\qquad\qquad\qquad\qquad\times K(\mu^{n-1},P(n-1))M(P(n-1),P)\Bigg)\end{align*}\begin{align*}
    &=\sum_{\substack{P(n-1),P(n-2):\\\sh P(n-1)=\mu^{n-1},\sh P(n-2)=\mu^{n-2}}}\Bigg(L(\emptyset,\mu^1)\dots K(\mu^{n-2},P(n-2))\\
    &\qquad\qquad\qquad\qquad\qquad\qquad\times M(P(n-2),P(n-1))M(P(n-1),P)\Bigg)\\
    &=\dots=\sum_{(P(i))_{i=1}^{n-1}:\ \sh P(i)=\mu^i}\Bigg(L(\emptyset,\mu^1)K(\mu^1,P(1))M(P(1),P(2))\times\\
    &\qquad\qquad\qquad\qquad\qquad\qquad\qquad\qquad\qquad\qquad\dots\times M(P(n-1),P)\Bigg).
  \end{align*}
  Now, from the 
  definition of $L$, $K$ and $M$, for $P(1)\in\vart_l$ that has only one entry $k$ and whose shape is $\mu^1=(1)$,
  \begin{align*}
    L(\emptyset,\mu^1)=1;\quad K(\mu^1,P(1))=M(\emptyset,P(1))=a_k.
  \end{align*}
This completes the proof.

\subsection{Proof of Proposition \ref{limpsi}}

Let $\lambda\vdash n$ and note that, for $l>n$, $\Delta_l(\lambda)=\Delta(\lambda)$.
We want to show that
$$\lim_{l\to\infty} \Psi_{(1/l)^l}(\lambda) = \frac{f^\lambda(q)}{n! (1-q)^n\Delta(\lambda)}.$$
From the definition of $\Psi_a$, this is equivalent to
$$\lim_{l\to\infty} l^{-n} \sum_{P\in\vart_l:\ \sh P=\lambda} \kappa(P)=
  \frac{f^\lambda(q)}{n! (1-q)^n\Delta(\lambda)}.$$
Write
$$ \sum_{P\in\vart_l:\ \sh P=\lambda} \kappa(P)=A+B$$
where $A$ denotes the sum over tableaux with distinct entries and $B$ denotes the remaining sum.
Assume $l>n$.  By \eqref{kapparho}, if $P$ has distinct entries, then
$$\kappa(P)=\frac{\rho(\hat P)}{(1-q)^n\Delta(\lambda)} .$$
Hence
$$l^{-n} A=l^{-n}\begin{pmatrix}l\\n\end{pmatrix}\sum_{Q\in\vars_n} \frac{\rho(Q)}{(1-q)^n\Delta(\lambda)} \to 
  \frac{f^\lambda(q)}{n! (1-q)^n\Delta(\lambda)}$$
as $l\to\infty$.  Thus it remains to show that $l^{-n} B\to 0$.  We first show that $\kappa(P)$ is bounded
for $P\in\vart_l$ with $\sh P=\lambda$.  To see this, observe that if $P$ has entries from the set $\{i_1,\ldots,i_m\}$
where $i_1<\cdots<i_m$ and $\tilde P$ denotes the tableau obtained from $P$ by replacing $i_k$ by $k$,
for each $k=1,\ldots,m$, then $\kappa(\tilde P)=\kappa(P)$.  It follows that
$$\kappa(P)\le \max_{T\in\vart_n}\kappa(T)<\infty.$$
Now, by the usual Robinson-Schensted correspondence, the number of $P\in\vart_l$ with $\sh P=\lambda$
which don't have distinct entries is at most the number of words $w\in [l]^n$ which don't have distinct entries,
and this is given by
$$N(l,n)=l^n-\begin{pmatrix}l\\n\end{pmatrix} n! .$$
Clearly, $l^{-n} N(l,n)\to 0$ as $l\to\infty$, so we are done.

\bigskip

{\em Acknowledgements.}  Research of N.O'C. supported in part by EPSRC grant number
EP/I014829/1. Research of Y.P. supported by EPSRC grant number EP/H023364/1.

\end{document}